\documentclass[12pt]{article}
\usepackage{verbatim}
\usepackage{graphicx}
\usepackage{hyperref,url,breakurl}
\usepackage{amsthm}
\usepackage{epsfig}
\title{Pan Galactic Division}
\author{Rich Schwartz}

\newtheorem{theorem}{Theorem}[section]

\begin{document}
\maketitle
\section{Introduction}

The purpose of these notes is to explain Peter Doyle and 
Cecil Qiu's
proof
\cite{doyleqiu:four}
that ``division by $N$'' is possible on a set-theoretic
level. To state the result formally, 
let $S \times N$ stand for $S \times \{0,...,(N-1)\}$.
Here $S$ is an arbitrary set.

\begin{theorem}
\label{doyle}
Let $A$ and $B$ be sets.
If there is an
injective map from $A \times N$ to $B \times N$
then there is an injective map from $A$ to $B$.
\end{theorem}

The significance of Doyle and Qiu's proof is that it
produces an effective (i.e. canonical) 
injection from $A \to B$ in a simple
and straightforward way.  In particular,
their proof avoids using the Well-Ordering Principle
or any other form of the Axiom of Choice.
For comparison,
in 1949 Tarski
\cite{tarski:cancellation}
gave an effective proof,
but Tarski's proof
is complicated.
Conway and Doyle
\cite{conwaydoyle:three}
sought a simpler proof,
but found nothing as simple as the proof here.
For more on the long and tangled
history of this problem,
see 
\cite{doyleqiu:four}.

Theorem \ref{doyle} is an immediate consequence of
induction and the following result.

\begin{theorem}
Let $N$ be any positive integer.  If there
is an injection from $A \times N$ to $A \times N$
then there is an injection from
$A \times (N-1)$ to $B \times (N-1)$.
\end{theorem}

The proof is the same regardless of the choice of $N$.
Following Doyle and Qiu, I'll give the proof when $N=4$, so that it
looks like a game of cards. Doyle and Qiu call their
game Pan Galactic Division, because they contend that a
definite fraction of intelligent civilizations in the
universe will have hit upon their canonical algorithm. 

\section{The Rules of the Game}

Think of $B \times 4$ as a deck of cards, where
the pictures on the cards are the members of $B$ and the
suits {\bf spades, hearts, diamonds, clubs\/}
are just other names for $3,2,1,0$.  
Think of $A$ as a set of players,
each staring
at $4$ spots in front of them, 
named {\bf spades, hearts, diamonds, clubs\/}
in order from left to right.

Enough 
of the cards in the deck are dealt out to the players
so that each player has a hand of $4$ cards, placed
into the $4$ spots.  (The suit on a card need not
match the suit which names the spot it is in.)
The injective map from
$A \times 4$ to $B \times 4$ just maps each pair
(player, spot) to the card in that location.
The game has alternating rounds, called {\it Shape Up\/} and 
{\it Ship Out\/}.
\newline
\newline
{\bf (1) Shape Up:\/}
Each player having at least one spade arranges to have
a spade as the leftmost card.  This is done by
swapping the leftmost spade with the leftmost card,
if such a swap is necessary. Figures 1 and 2 show
an example, before and after the player shapes up.

\begin{center}
\resizebox{!}{1.5in}{\includegraphics{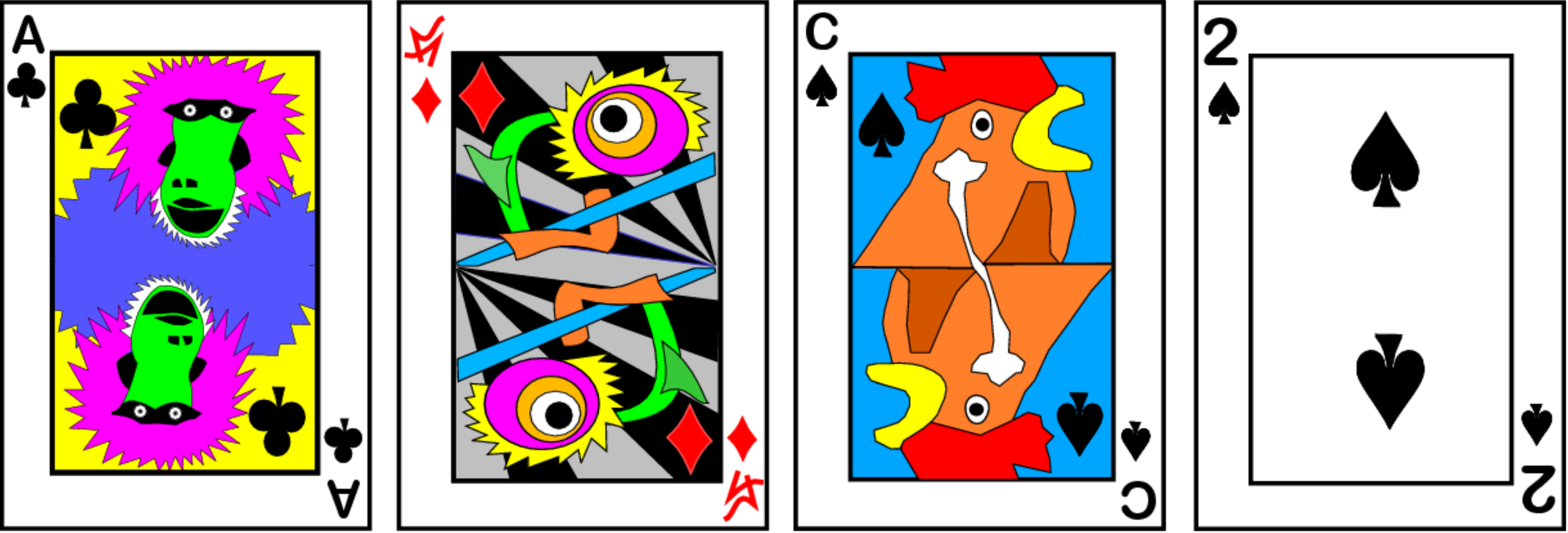}}
\newline
Figure 1: Before shaping up: The Chicken of spades is the leftmost spade.
\newline
\end{center}

\begin{center}
\resizebox{!}{1.5in}{\includegraphics{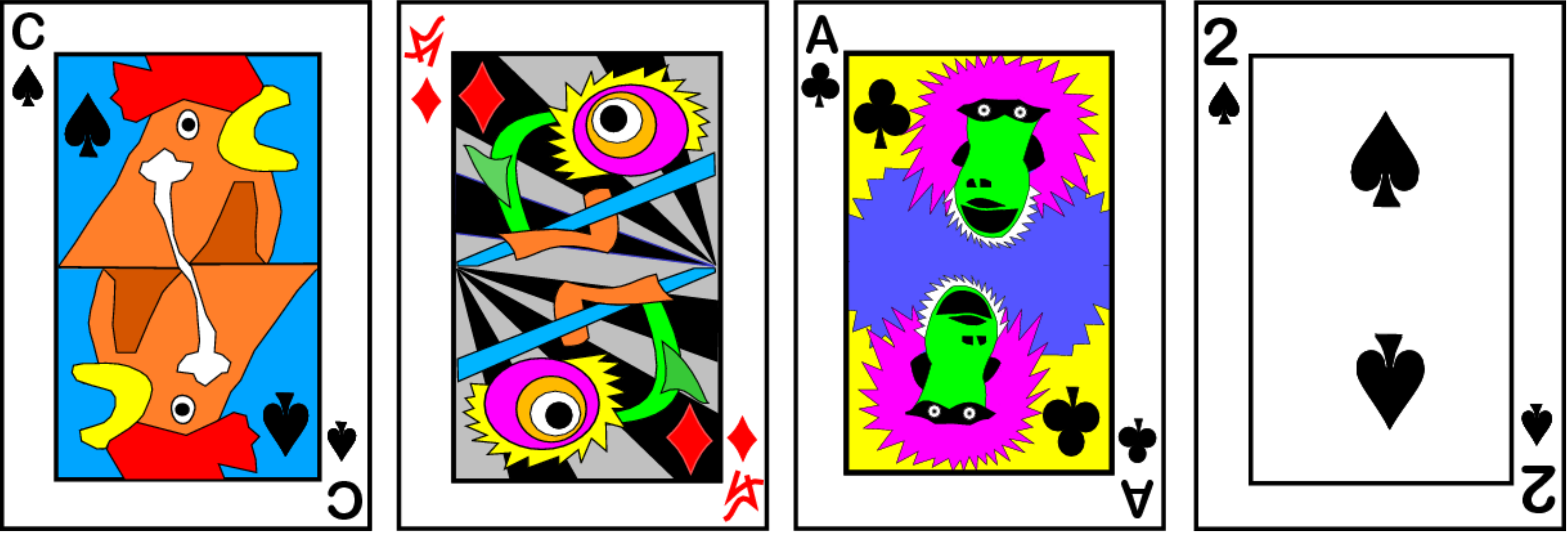}}
\newline
Figure 2: After shaping up:  The Chicken of spades swaps with Ape of clubs.
\newline
\end{center}

\noindent
{\bf Names of the Hands:\/}
All players shape up simultaneously.
Once that is done, any hand which has a
spade in it will also have a spade in the leftmost
spot.  In this case, the name of this leftmost
spade names the hand.  In our example, the name of
the hand is ``Chicken'', because the Chicken of
spades occupies the leftmost position.
(Since the sets $A$ and $B$ might be very large,
we can't expect our cards to be as in a traditional deck.)
The hands without spades
are not named.
\newline
\newline
{\bf (2) Ship Out\/}:
After a player shapes up, any named hand
in her hand
will have a spade on the left.
Say that a {\it bad spade\/} is a spade in a named
hand which is not the leftmost spade.  In our
example above, the $2$ of spades is a bad spade.
A hand can have up to three bad spades.

In the Ship Out round, each player with a bad spade
swaps their leftmost bad spade with another
card in the game, according to the following
rule.  Recall that the spots in front of the
player are named spades, hearts, diamonds, clubs,
from left to right.
The player asks for the card which has the same
suit as the name of the spot occupied by the
leftmost bad spade, and the same name
as the name of the hand.  There is never
any conflict: two players will not ask
for the same card.

Let's continue with our example above.
In our example, there is one bad spade, namely
the $2$ of spades.  Even though this card 
lies in the rightmost spot, it is the leftmost
bad spade, because it is the only bad spade.
The $2$ of spades occupies the clubs spot,
according to our scheme, and it lies in the
Chicken hand.  Therefore, the player swaps
the $2$ of spades with the Chicken of clubs.

The Chicken of clubs either is part of
another player's hand, or elsewhere in the player's own hand,
or else is still in the
deck; the distinction does not matter.  
Figure 3 shows our hand after the player ships out.

\begin{center}
\resizebox{!}{1.5in}{\includegraphics{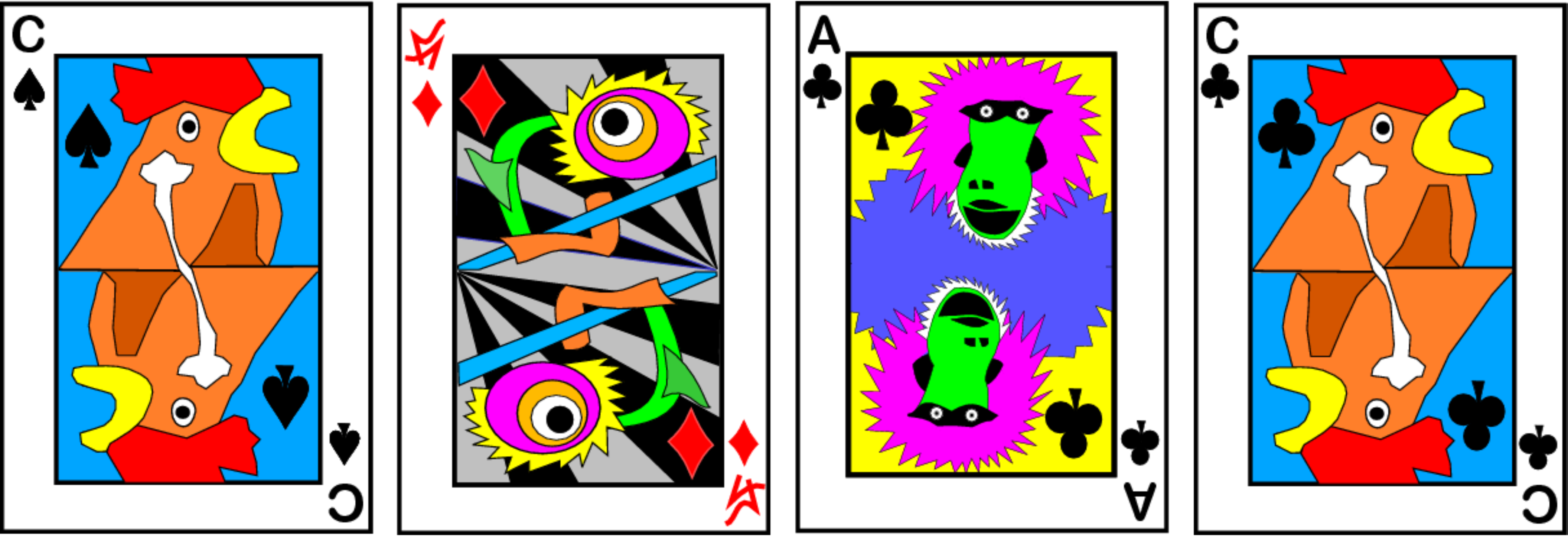}}
\newline
Figure 3: Two of spades is traded for the Chicken of clubs.
\newline
\end{center}

\section{Playing the Game}

The game is now played 
indefinitely, with the rounds alternating: 
shape up, ship out, shape up, ship out, \ldots
Notice that a
hand never loses its leftmost spade card,
because 
the Ship Out rule will never require anyone to call for a spade.
So, any player shapes up at most once.

The quality of a named hand is assessed according to
how many cards have names which match the name of the
hand and suits which match the name of the spot they
occupy. Let's consider our example. Before
shipping out, our player has $1$-of-a-kind, so to
speak, because the Chicken of spades is her only card
which has the same name as the hand (by definition) and
is in the correct position.  After shipping out,
our player has acquired the Chicken of clubs and placed
it into the clubs position. Thus, now she has $2$-of-a-kind.
Notice that our player has improved her hand by 
shipping out.

As our example suggests, a player always improves her
hand by shipping out a spade.   A player's hand can
have a spade {\it shipped in\/} by having one of 
her cards called away.
For instance, in the example below, our player
loses the ape of clubs and acquires the bolt
of spades.  The bolt of spade is a bad spade in the
diamonds position, as shown in Figure 4.

\begin{center}
\resizebox{!}{1.7in}{\includegraphics{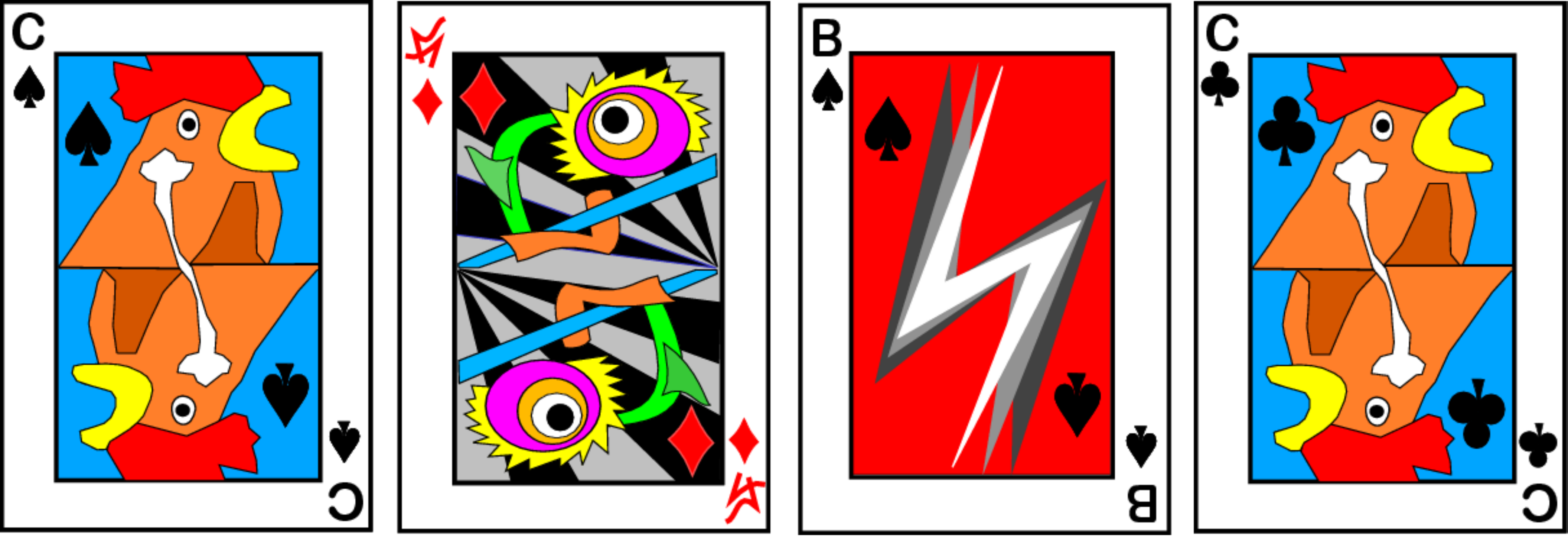}}
\newline
Figure 4: The bolt of spades occupies the diamonds position.
\end{center}

In the next round, our player ships out
the bolt of spades and gets the Chicken of diamonds,
as shown in Figure 5.
In this way, she has improved her hand to $3$-of-a-kind.
So, any time a player participates in a Shipping Out round,
either actively or passively, her hand improves 
after at most two more rounds of the game.

\begin{center}
\resizebox{!}{1.7in}{\includegraphics{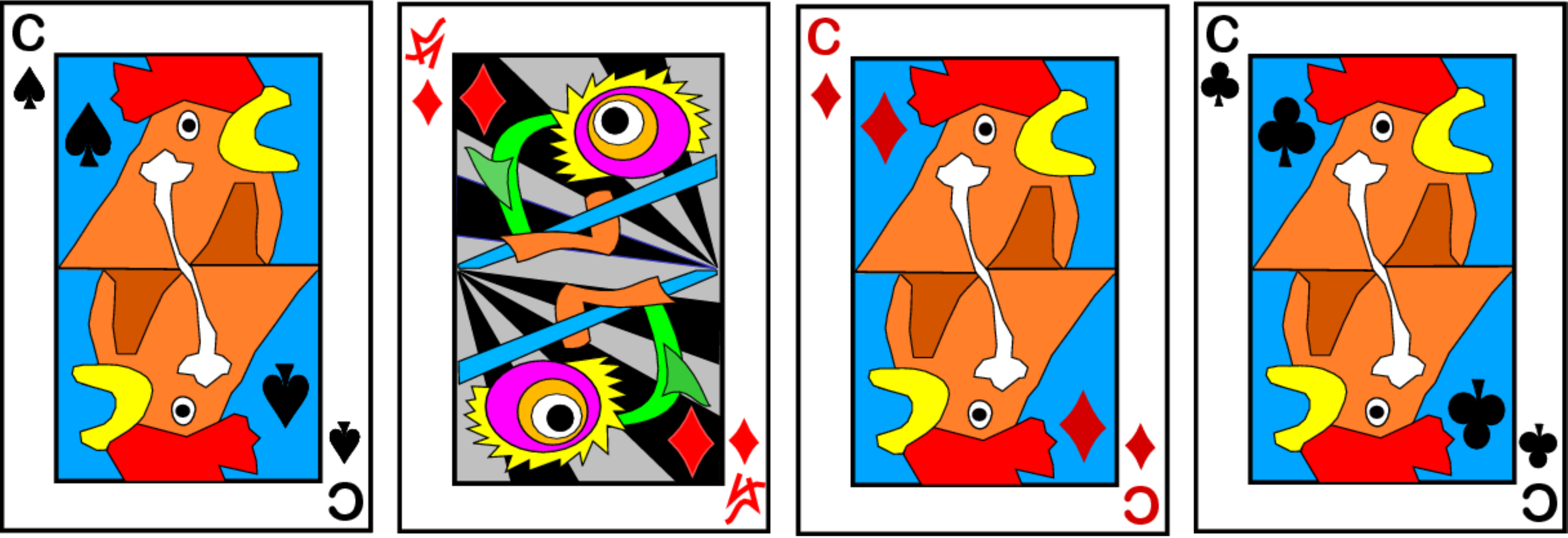}}
\newline
Figure 5: The bolt of spades is traded for the Chicken of diamonds.
\end{center}

\section{Conclusion}

From the discussion above, we conclude that during this
indefinite game, a player's hand
can change at most $8$ times.  Moreover, once a player's
hand has stopped changing, the player has no bad spades.
In other words, once the hand has stabilized, it has
no spades in the $3$ non-left spots.

We imagine that the game proceeds indefinitely.
There is a well-defined non-spade card associated
to each pair (player, non-left spot):  We simply
wait until the player's hand stabilizes - and it
will stabilize - and then we see what card is
occupying the relevant spot.  This gives us
an injective map from $A \times 3$ to $B \times 3$.

\bibliography{four}
\bibliographystyle{hplain}
\end{document}